\documentclass[11pt,a4paper
]{article}
\usepackage{amsthm,amsfonts,amsmath,amscd,amssymb}
\usepackage{latexsym}
\usepackage{euscript}
\usepackage{enumitem}
\usepackage{graphicx}
\usepackage{tikz}
\usepackage{caption}
\usepackage{subcaption}
\usepackage{cite}
\usepackage[utf8]{inputenc}
\usepackage[english]{babel}
\usepackage{xcolor}

\usepackage{hyperref}
\usepackage{etoolbox}
\usepackage{forloop}

\hypersetup{colorlinks=true}
\numberwithin{equation}{section}
\numberwithin{figure}{section}
\numberwithin{table}{section}

\newcommand{\key}[1]{\par\vskip 0.5em\small{\slshape Key words:\/} #1}

\newcommand{\msc}[1]{\par\vskip 0.5em\small{\slshape Mathematical Subject Classification 2020:\/} #1}

\newcommand\arXiv[1]{\url{https://arxiv.org/abs/#1}}

\theoremstyle{plain}
\newtheorem{theorem}{Theorem}

\newtheorem{lemma}[theorem]{Lemma}
\theoremstyle{definition}
\newtheorem{definition}[theorem]{Definition}
\theoremstyle{remark}

\newtheorem{example}[theorem]{Example}

\usepackage{bm}

\def\fracd{\displaystyle\frac}
\def\sumd{\displaystyle\sum}

\def\prodd{\displaystyle\prod}
\def\capd{\displaystyle\bigcap}
\def\cupd{\displaystyle\bigcup}

\title{EIGENVALUE DISTRIBUTION OF
LARGE WEIGHTED RANDOM SPARSE UNIFORM $q$-HYPERGRAPHS}

\author{Valentin Vengerovsky\footnote{B. Verkin Institute for Low Temperature Physics and Engineering of the National Academy of Sciences of Ukraine, 47 Nauky Ave., Kharkiv, 61103, Ukraine. E-mail: \href{mailto:vengerovsky@ilt.kharkov.ua}{vengerovsky@ilt.kharkov.ua}}}

\begin{document}
\maketitle



\newcommand{\Compl}{\mathbb{C}}
\newcommand{\R}{\mathbb{R}}
\newcommand{\Z}{\mathbb{Z}}

\newcommand{\conseq}{\Rightarrow}

\newcommand{\Mat}{\mathrm{Mat}}
\newcommand{\conj}[1]{\overline{#1}}
\newcommand{\transp}[1]{#1^T}
\newcommand{\rev}[1]{\dual{#1}}
\newcommand{\dual}[1]{#1^R}
\newcommand{\cConjScl}[1]{\bar{#1}}
\newcommand{\cConjMat}[1]{#1^*}
\newcommand{\aConjScl}[1]{#1^*}
\newcommand{\aConjMat}[1]{#1^+}
\newcommand{\ds}[1]{\check{#1}}

\newcommand{\setcharf}[1]{\mathbbm{1}_{#1}}
\newcommand{\charfp}{\psi}
\newcommand{\asKer}[2]{K\left(#1, #2\right)}
\newcommand{\CF}{\mathsf{f}}

\newcommand{\Gin}{\mathrm{Gin}}
\newcommand{\ens}{M_n}
\newcommand{\sclA}{y}
\newcommand{\matA}{Y}

\newcommand{\cSclGin}{x}
\newcommand{\cMatGin}{X}
\newcommand{\cMatPos}{\mathcal{Z}}
\newcommand{\cSclN}{t}
\newcommand{\cVecN}{\bm{\cSclN}}
\newcommand{\cSclA}{q}
\newcommand{\cMatA}{Q}
\newcommand{\cSetMatA}{\bm{Q}}
\newcommand{\cVecGI}{\bm{h}}
\newcommand{\cSclB}{a}

\newcommand{\aSclA}{\xi}
\newcommand{\aVecA}{\bm{\aSclA}}
\newcommand{\aMatA}{\Xi}
\newcommand{\aSetMatA}{\mathbf{\aMatA}}
\newcommand{\aSclB}{\phi}
\newcommand{\aVecB}{\bm{\aSclB}}
\newcommand{\aMatB}{\Phi}
\newcommand{\aSclBt}{\varphi}
\newcommand{\aVecBt}{\aVecB}
\newcommand{\aSclC}{\theta}
\newcommand{\aVecC}{\bm{\aSclC}}
\newcommand{\aMatC}{\Theta}
\newcommand{\aSclCt}{\vartheta}
\newcommand{\aVecCt}{\aVecC}
\newcommand{\aSclD}{\rho}
\newcommand{\aVecD}{\bm{\aSclD}}
\newcommand{\aSclE}{\tau}
\newcommand{\aSclF}{\nu}
\newcommand{\aVecF}{\bm{\aSclF}}
\newcommand{\aSclG}{\upsilon}
\newcommand{\aVecG}{\bm{\aSclG}}
\newcommand{\aMatG}{\Upsilon}
\newcommand{\aSclH}{\aSclE}
\newcommand{\aVecH}{\bm{\aSclH}}
\newcommand{\aVecGrI}{\bm{\upsilon}}

\newcommand{\tempcumul}{\kappa}
\newcommand{\cumul}[2]{\tempcumul_{#1,#2}}
\newcommand{\realcumul}{\tempcumul}
\newcommand{\seccumul}{\cumul{2}{0}}
\newcommand{\indexset}{\mathcal{I}}
\newcommand{\emptyindex}{\varnothing}
\newcommand{\diLambda}{\mathcal{M}}
\newcommand{\stpointsnbh}{\Omega_n}
\newcommand{\idom}{\mathcal{D}}
\newcommand{\Vanddet}{\triangle}
\newcommand{\herm}{\mathcal{H}}
\newcommand{\USp}{\mathrm{USp}}
\newcommand{\scoeff}{\mathsf{c}}
\newcommand{\vol}{\mathrm{vol}}
\newcommand{\skB}{\mathcal{B}}
\newcommand{\Hess}{\mathsf{H}}
\newcommand{\Fperm}{\mathcal{F}}

\newcommand{\partder}[2]{\frac{\partial #1}{\partial #2}}
\newcommand{\der}[2]{\frac{d #1}{d #2}}

\newcommand{\abs}[1]{\left\lvert#1\right\rvert}
\newcommand{\abssized}[2][ ]{#1\lvert#2#1\rvert}
\newcommand{\norm}[1]{\left\lVert#1\right\rVert}
\newcommand{\normsized}[2][ ]{#1\lVert#2#1\rVert}

\newcommand{\tr}{\mathop{\mathrm{tr}}}
\newcommand{\Pf}{\mathop{\mathrm{Pf}}}
\newcommand{\E}{\operatorname{\mathbf{E}}}
\newcommand{\diag}{\mathop{\mathrm{diag}}}



\begin{abstract}
We study eigenvalue distribution of the adjacency
matrix $A^{(N,p,q)}$
of weighted random uniform $q$-hypergraphs $\Gamma= \Gamma_{N,p,q}$. We assume
that the graphs have
$N$ vertices and the average number of hyperedges attached to one vertex  is
$(q-1)!\cdot p$. To each edge of the graph $e_{ij}$ we assign a weight given by a random variable
$a_{ij}$ with all moments finite. We consider the moments of
normalized eigenvalue counting function $\sigma_{N,p,q}$  of
$A^{(N,p,q)}$. Assuming all moments of $a$ finite, we
 obtain recurrent relations that determine the
moments of the limiting measure $\sigma_{p,q} = \lim_{N\to\infty}
\sigma_{N,p,q}$.

\key{ sparse random uniform $q$-hypergraph, normalized eigenvalue counting measure, limiting mesure.}

\msc{60B20, 15B52.}
\end{abstract}


\section{Introduction}

In recent years, interest in the spectral properties of ensembles of sparse random matrices has sharply increased. It is anticipated that the spectral characteristics of sparse random matrices will differ from those of ensembles of most matrices with independent elements (see \cite{W}, along with the survey works \cite{KKPS}, \cite{Pa:00}, and the literature cited therein). The investigation of various issues for different ensembles of sparse random matrices can be seen in the works \cite{BG1}, \cite{BG2}, \cite{B}, \cite{KV}, \cite{KSV}, \cite{PJ}, \cite{CS}, \cite{V0}, \cite{V},
\cite{V1}, \cite{V2}, \cite{KKP}, \cite{APS2}.

Spectral properties of different hypergraphs have also been actively studied in the last decades \cite{DZ}, \cite{SZ},\cite{SB}.

In this paper, we investigate a limiting normalized eigenvalue counting  measure of large sparse random weighted uniform $q$-graphs.


\section{Main results}


    We can introduce the randomly weighted
adjacency matrix of random sparse uniform $q$-hypergraphs.
Let $\Xi=\{a_{e} ,\; e \subset \!{\mathbb N}, |e|=q\}$ be
the  set of
jointly independent identically distributed (i.i.d.) random
variables determined
  on the same probability space and possessing the moments
\begin{equation}\label{const_all_mom}
  {\mathbb E}a^{k}_{e}\!=\!X_{k}\!<\!\infty \qquad
  \forall \; e \in {{\mathbb N}\choose q},
\end{equation}
where ${\mathbb E}$ denotes the mathematical expectation corresponding to
$\Xi$.

Given $ 0\!<p\!\leq \!N$, let us define the family
  $D^{(p,q)}_N\!=\!\{d^{(N,p,q)}_{e},
\; e\in {{ N}\choose q}\}$ of jointly independent
random variables
\begin{equation}
d^{(N,p,q)}_{e}\!=\! \left\{ \begin{array}{ll} 1,&
\textrm{with} \ \textrm{probability } \ p/N^{q-1} ,
\\0,& \textrm{with} \ \textrm{probability} \ 1-p/N^{q-1} ,\\ \end{array}
\right.
\end{equation}
We  assume that $ D^{(p,q)}_N$
is independent from $\Xi$.

  Now we can consider the real symmetric $N\times N$ matrix
$A^{(N,p,q)}(\omega)$:
\begin{equation}\label{dilute}
\left[A^{(N,p,q)}\right]_{ij}\!=\!
\sumd_{ e\in {{ N}\choose q}} A^{(e,N,p,q)}_{ij},
\end{equation}
where
\begin{equation}
 A^{(e,N,p,q)}_{ij}\!=\! \left\{ \begin{array}{ll} a_e\cdot d^{(N,p,q)}_e,&
\textrm{if} \ i\in e, j\in e, i\neq j  ,
\\0,& \textrm{otherwise.} \\ \end{array}
\right.
\end{equation}
Let weighted adjacency matrix $A^{(N,p,q)}$  has $N$ real eigenvalues
$\lambda^{(N,p,q)}_1\!\leq\!\lambda^{(N,p,q)}_2 \!\leq\!\ \ldots
\!\leq\!\ \lambda^{(N,p,q)}_N$.

The normalized eigenvalue counting function of  $A^{(N,p,q)}$ is determined by the formula
 $$
\sigma\left({\lambda; A^{(N,p,q)}}\right)\!=\!\frac{\#
\left\{j:\lambda^{(N,p,q)}_j\!<\!\lambda\right\}}{N}.
$$

\begin{theorem}\label{main_thm}
Under condition
\begin{equation}\label{Karleman_cond}
X_{2m}\le \left(C \cdot m\right)^{2m}, m\in \mathbb{N}
\end{equation}
 measure $\sigma\left({\lambda; A^{(N,p,q)}}\right)$ weakly converges in probability to nonrandom measure  $\sigma_{p,q}$
\begin{equation}
\sigma\left({\cdot \ ; A^{(N,p, q)}}\right)\to \sigma_{p,q},\ N \to\infty ,
\end{equation}
which can be uniquely determined  by its moments
\begin{equation}\label{5}
\int \lambda^k d \sigma_{p,q}=
m^{(p,q)}_{k}\!=\!\sum^{k/2}_{i=0} S(k,i),
\end{equation}
where numbers  $S(k,i)$, can be found from the following system of recurrent relations
\begin{multline}\label{ms_r}
  \mathrm{S}(l,r)=p\displaystyle\sum_{F=1}^{l} X_F \cdot  \sum_{\varkappa=1}^{q} \sum_{\substack{f\in {\mathbb Z}_{+}^\varkappa \\ ||f||_1=F \\f_1\le r}}\frac{1}{(q-\varkappa)!}\cdot K^{(\varkappa,q)}_{1;f}\cdot
         \\
          \sum_{u=0}^{l-f}{r-1\choose f_1-1} \mathrm{S}(u_1,r-f_1)
       \cdot \prod_{i=2}^\varkappa \sum_{u_i=0}^{l-F-u_1-\ldots -u_{i-1}}\sum_{v_i=0}^{u_i/2} {f_i+v_i-1\choose f_i-1} \cdot
       \mathrm{S}(u_i,v_i),
\end{multline}

with following initial conditions
\begin{equation}\label{ini_cond}
S(l,0)= \delta_{l,0}.
\end{equation}
$\{K^{(\varkappa,q)}_{j;f}\}$ yield the next system of recurrent relations
\begin{equation}\label{sys_K}
K^{(\varkappa,q)}_{j;f}=\sumd_{\substack {a=1,a\neq j\\ f-1_a\in\mathbb{Z}_+^q}}^{\varkappa}K^{(\varkappa,q)}_{a;f-1_a}+\delta_{\varkappa -j}\cdot \sumd_{\substack{a=1 \\f-1_a\in\mathbb{Z}_+^q}}^{\varkappa -1}K^{(\varkappa -1,q)}_{a;f-1_a}, j\in[\varkappa]
\end{equation}
with following initial conditions
\begin{equation}\label{ini_cond_K}
K^{(1,q)}_{j;f}=\delta_{j-1}\prod_l^q \delta_{f_l}.
\end{equation}
\end{theorem}

The following denotations are used:
$$
{\mathcal{M}}^{(N,p,q)}_t\!=\! \int \lambda^t d \sigma\left({\lambda;
A^{(N,p,q)}}\right), \ M^{(N,p,q)}_t=\mathbb{E}{\cal{M}}^{(N,p,q)}_t,
$$
$$
C^{(N,p,q)}_{t,m}= \mathbb{E}\left\{{\mathcal{M}}^{(N,p,q)}_{t}\cdot {\cal{M}}^{(N,p,q)}_{m} \right\}-\mathbb{E}\left\{{\cal{M}}^{(N,p,q)}_{t} \right\}\cdot \mathbb{E}\left\{{\cal{M}}^{(N,p,q)}_{m} \right\} .
$$
Theorem \ref{main_thm} is a corollary of Theorem \ref{thm:1}
\begin{theorem} \label{thm:1}
Assuming  conditions (\ref{Karleman_cond}),
(i) Correlators $C^{(N,p,q)}_{k,m}$  vanish in the limit $\ N\to \infty$:
\begin{equation}
C^{(N,p,q)}_{k,m}\le {\displaystyle\frac{C(k,m,p,q)}{N}} ,  \ \forall\  k,m \in \mathbb{N}.
\end{equation}

(ii) The limit of s-th moment exists for all $s\in \mathbb{N}$:
\begin{equation}\label{5}
  \lim_{N \to \infty}M^{(N,p, q)}_k =
m^{(p,q)}_{k}\!=\!\sum^{k/2}_{i=0} S(k,i),
\end{equation}
where  numbers
  $S(k,i)$
are determined by  (\ref{ms_r}) -(\ref{ini_cond_K}).

(iii) The limiting moments $\left\{m^{(p,q)}_{k}\right\}_{k=1}^{\infty}$ obey Carleman's condition
\begin{equation}\label{Carleman}
{\displaystyle\sum_{k=1}^{\infty}}{\displaystyle\frac{1}{\sqrt[2k]{m^{(p,q)}_{2k}}}}=\infty
\end{equation}
\end{theorem}

\section{Proof of Theorem 1}
\subsection{Walks and contributions}

\quad\ Using independence of families $\Xi$ and
$D^{(p,q)}_N$, we have {\setlength\arraycolsep{1pt}
\begin{eqnarray}
  M^{(N,p,q)}_k\!&=&\!\int {\mathbb E} \{ \lambda^k d \sigma_{A^{(N,p,q)}}
\} \!=\!{\mathbb E} \left(\frac{1}{N}\sum_{i=1}^N
[\lambda^{(N,p, q)}_i]^k \right)\!=\! \frac{1}{N} {\mathbb E}\left(Tr
[A^{(N,p,q)}]^k\right) \!=\! \nonumber \\ & =&\! \frac{1}{N}
\sum^{N}_{j_1=1} \sum^{N}_{j_2=1} \ldots
  \sum^{N}_{j_{k}=1} {\mathbb E} \left( A^{(N,p,q)}_{j_1,j_2}
  A^{(N,p,q)}_{j_2,j_3} \ldots A^{(N,p,q)}_{j_{k},j_1}
  \right) \!=\! \nonumber \\
\label{base} &=&\! \frac{1}{N}
 \sum_{\substack{j_1, j_2, \ldots, j_k\in [n] \\ e_1,e_2,\ldots, e_k\in {[n]\choose q} \\ j_1,j_2\in e_1,\ldots, j_k,j_1\in e_k}}
{\mathbb E} \left( A^{(N,p,q)}_{j_1,j_2;e_1}
  A^{(N,p,q)}_{j_2,j_3;e_2} \ldots A^{(N,p,q)}_{j_{k},j_1;e_k}
  \right) =
\nonumber \\ & &
 = \frac{1}{N}
 \sum_{\substack{j_1, j_2, \ldots, j_k\in [n] \\ e_1,e_2,\ldots, e_k\in {[n]\choose q} \\ j_1,j_2\in e_1,\ldots, j_k,j_1\in e_k}}
{\mathbb E} \left( a_{e_1}\cdot a_{e_2}\cdot\ldots \cdot a_{e_k}\right)\cdot
\nonumber \\ & &
  \cdot{\mathbb E} \left( d^{(N,p,q)}_{e_2}\cdot d^{(N,p,q)}_{e_1}\cdot\ldots\cdot d^{(N,p,q)}_{e_k}\right).
 \end{eqnarray}

Let $W^{(N)}_{k}$ be a set of closed walks of $k$ steps over
the set $[N]$:
$$
W^{(N)}_{k}\!=\!\{w\!=\!(w_1,e_1,w_2,e_2,\ldots,w_k,e_k,w_{k+1}=w_1):
\forall i \!\in\! [k] \;\: w_i\!\in\!
[N], \;\: e_i\!\in\!
{ [N]\choose q}, \;\: j_i,j_{i+1}\!\in\!
e_i\}.
$$
For $w\!\in\!W^{(N)}_k$ let us denote
  $a(w)\!=\!\prodd_{i=1}^{k} a_{e_i}$,
 $d^{(N,p,q)}(w)\!=\!\prodd_{i=1}^{k} d^{(N,p,q)}_{e_i}$.
Then we have
\begin{equation}\label{m_ms1}
M^{(N,p,q)}_k\!=\!\frac{1}{N} \sum_{w\in W^{(N)}_k} {\mathbb E} a(w)
\cdot {\mathbb E} d^{(N,p,q)}(w).
\end{equation}
Let $w\!\in\! W^{(N)}_k$
  and
$e \!\in\! {[N]\choose q}\ $. Denote by $n_w(e)$ the number
of steps along hyperedge $e$ in closed walk $w$
$$
n_w(e)=\#\{i \!\in\! [k]:\; e_i\!=\! e\}.
$$
Then
$$
{\mathbb E}a(w)\!=\! \ \prod_{e} X_{n_w(e)}.
$$

  Given $w\!\in\! W^{(N)}_k$,
let us define the sets $V_w=\cup_{i=1}^{k}e_k$ and
$E_w=\cup_{i=1}^{k}\{e_k)\}$. It is easy to see that $G_w\!=\!(V_w,E_w)$ is a
simple uniform $q$-hypergraph and the walk $w$ covers the hypergraph
$G_w$. Let us call $G_w$ the skeleton of walk $w$. Then we obtain
$$
  {\mathbb
E}a(w)\!=\!\prod_{e\in E_w} {\mathbb E}a^{n_w(e)}_e\!=\! \prod_{e\in
E_w} X_{n_w(e)}.
  $$
Similarly we can write
$$
  {\mathbb E}d^{(N,p,q)}(w)\!=\!\prod_{e\in E_w} {\mathbb
E}\left([d^{(N,p,q)}_e]^{n_w(e)} \right)\!=\! \prod_{e\in E_w}
\frac{p}{N^{q-1}}.
$$
Then, we can rewrite (\ref{m_ms1}) in the form
$$
M^{(N,p,q)}_k\!=\!\frac{1}{N}\sum_{w\in W^{(N)}_k}  \prod_{e\in E_w}
\frac{p\cdot X_{n_w(e)}}{N^{q-1}}\!=
$$
\begin{equation}\label{m_ms2}
=\!\sum_{w\in W^{(N)}_k} \left(\frac{p^{|E_w|}}{N^{(q-1)\cdot|E_w|+1}
}\prod_{e\in E_w}X_{n_w(e)} \right)\!=\!\sum_{w\in
W^{(N)}_k}\theta(w),
\end{equation}
where $\theta(w)$ is the contribution of the walk $w$ to the
mathematical
  expectation of the corresponding moment. To perform the limiting
transition
   $N\to\infty$ it is natural to separate $W^{(N)}_k$ into classes of
equivalence.
   Walks $w^{(1)}$ and $w^{(2)}$ are equivalent
  $ w^{(1)}\sim w^{(2)},\;$
if and only if there exists a bijection $\phi$ between the sets of
vertices $V_{w^{(1)}}$ and  $V_{w^{(2)}}$ such that for
$i\in [k]
  \;\; w^{(2)}_i\!\!=\!\!\phi(w^{(1)}_i), e^{(2)}_i\!\!=\!\!\phi(e^{(1)}_i)$.
$ w^{(1)}\sim w^{(2)}\;\Longleftrightarrow
\; \exists \phi: \ V_{w^{(1)}}\stackrel{bij}{\to}
  V_{w^{(2)}}:
  \;\forall \;i \!\in\! [k];
\; w^{(2)}_i\!\!=\!\!\phi(w^{(1)}_i), e^{(2)}_i\!\!=\!\!\phi(e^{(1)}_i)$.
   Let us denote by $[w]$ the class of equivalence of walk $w$ and by
$C^{(N)}_k$ the set of such classes. It is obvious that if two
walks $w^{(1)}$ and $w^{(2)}$ are equivalent then their
contributions are equal:
$$
  w^{(1)}\sim w^{(2)}\;\Longrightarrow
\theta(w^{(1)})\!=\!\theta(w^{(2)}).
$$
Cardinality of the class of equivalence $[w]$ is equal the number
of all mappings $\phi:V_w \to  [N]$ divided by the number of mappings that leave the closed walk unchanged. So it  is equal to the number
$\prodd_{i=1}^{|V_w|} N \cdot (N-1) \cdot \ldots \cdot (N-|V_{w}|+1)\cdot\prodd_{\substack{\alpha\subset 2^{[E_w]}\\ \alpha \neq \emptyset}}\frac{1}{|\beta_\alpha|!}$, where $\beta_\alpha=\left(\capd_{i\in\alpha}e_i \setminus \cupd_{i\in\overline{\alpha}}e_i\right)\setminus \widetilde{V}_w$ and $\widetilde{V}_w= \cupd_{i\in [k]}w_i$. Then we can rewrite
(\ref{m_ms2}) in the form
$$M^{(N,p)}_k\!=\!\sum_{w\in W^{(N)}_k}
\left(\frac{p^{|E_w|}}{N^{(q-1)\cdot|E_w|+1} }\prod_{e\in
E_w}X_{n_w(e)} \right)\!=
$$
$$
=\!\sum_{[w]\in C_k}
\frac{p^{|E_w|}}{N^{(q-1)\cdot|E_w|+1} }\prod_{e\in
E_w}X_{n_w(e)} \!\cdot
$$
\begin{equation}\label{m_ms3}
\left. \cdot\prodd_{i=1}^{|V_w|} N \cdot (N-1) \cdot \ldots \cdot (N-|V_{w}|+1)\cdot\prodd_{\substack{\alpha\subset 2^{[E_w]}\\ \alpha \neq \emptyset}}\frac{1}{|\beta_\alpha(w)|!}\right)\!=\!\sum_{[w]\in
C^{(N)}_k} \hat{\theta}([w]).
\end{equation}
In the second line of (\ref{m_ms3}) for every class $[w]$ we choose arbitrary walk $w$ corresponding to this class of equivalence.
\subsection{Minimal and essential walks}

   Class of walks $[w]$ of $C^{(N)}_k$ has at most k vertices.
    Hence, $C^{(1)}_k
\subset C^{(2)}_k \subset \ldots \subset C^{(i)}_k \subset
\ldots C^{(k_q)}_k = C^{(k_q+1)}_k= \ldots$. It is natural to denote
  $C_k=C^{(k_q)}_k$. Then (\ref{m_ms3}) can be written as
\begin{equation}\label{m_ms5}
m^{(q)}_k \!=\!\lim_{N \to \infty} \sum_{[w]\in C_k}
\left(N^{|V_w|-(q-1)\cdot|E_w|-1}\prod_{e\in
E_w}\fracd{ X_{n_w(e)}}{p^{-1}}\cdot\prodd_{\substack{\alpha\subset 2^{[k]}\\ \alpha \neq \emptyset}}\frac{1}{|\beta_\alpha(w)|!}\right).
\end{equation}
The set $C_k$ is finite. Regarding this and (\ref{m_ms5}), we
conclude that the class $[w]$ has non-vanishing contribution,
only if $|V_w|-(q-1)\cdot|E_w|-1 \!\geq \! 0$. But for each simple connected
uniform $q$-hypergraph $G=(V,E)$ the following inequality holds  $|V|\! \leq \! (q-1)\cdot|E|+1$, and the equality takes
place if and only if the hypergraph $G$ is
  a hypertree.

It is convenient to deal with $\widetilde{W}^{(N)}_{k}$ instead of $W^{(N)}_{k}$, where $\widetilde{W}^{(N)}_{k}$
 is a set of closed walks over the set  $[N] \cup \{\widetilde{1},\widetilde{2},\ldots, \widetilde{N}\}$.
 Let us consider $\widetilde{C}^{(N)}_{k}$, the set of equivalence  classes of $\widetilde{W}^{(N)}_{k}$. As a  representative of the equivalence  class $[w]\in \widetilde{C}^{(N,)}_{k}$, we can take a minimal  walk.
\begin{definition}
   A  closed walk  $w\in \widetilde{W}^{(N)}_{k}$ is called minimal if and only if f the following conditions are met. All vertices of $\widetilde{V}_w$ belong $\mathbb{N}$, and all vertices of $V_w\setminus \widetilde{V}_w$ belong the set $\{\widetilde{1},\widetilde{2},\widetilde{3}, \dots\}$. At each stage of the passage a new vertex is the minimum element among the unused numbers of $\mathbb{N}$. The set $V_w\setminus \widetilde{V}_w$ is is naturally partitioned into disjoint sets $\displaystyle \bigcup_{\substack{\alpha\subset 2^{[E_w]}\\ \alpha \neq \emptyset}} \beta_\alpha$. The sets $\displaystyle \bigcup_{\substack{\alpha\subset 2^{[E_w]}\\ \alpha \neq \emptyset}} \beta_\alpha$ are taken in lexicographic order and filled without repetitions in ascending order with elements of the set of natural numbers with a tilde.
\end{definition}
Let us denote the set of all minimal walks of $\widetilde{W}^{(N)}_{k}$ by
$\mathfrak{M}^{(N)}_{k}$.

\begin{example}\label{min_example} The sequence  $(1,\{1,2,4\},2,\{2,3,4\},3,\{2,3,4\},4,\{4,5,\widetilde{1}\},5,\{4,5,\widetilde{1}\},4,\{1,2,4\},1,$ $\{1,6,\widetilde{1}\},6,\{6,7,\widetilde{2}\},7,\{6,7,\widetilde{2}\},6,$ $\{1,6,\widetilde{1}\},1,
\{1,2,4\},4,\{2,3,4\},2,\{2,3,4\},3,$ $\{2,8,\widetilde{2}\},8,$ $\{8,9,\widetilde{1}\},9,$ $\{8,9,\widetilde{1}\},8,\{3,8,\widetilde{2}\},3,\{2,3,4\},2,\{1,2,4\},1,\{1,10,\widetilde{3}\},10,\{1,10,\widetilde{3}\},1)$  is a   minimal walk.\end{example}

\begin{definition}
The  minimal walk $w$ that has a hypertree as
a skeleton is an essential walk.
\end{definition}
The walk from example \ref{min_example} is not essential. Let us denote the set of all essential  walks of $\widetilde{W}^{(N)}_{k}$ by
$\mathfrak{E}^{(N)}_{k}$.
 Therefore we can rewrite
(\ref{m_ms5}) in the form
\begin{equation}\label{m_ms6}
m^{(p,q)}_k \!=\! \sum_{w\in \mathfrak{E}_k}
\prod_{e\in E_w}\frac{p\cdot  X_{n_w(e)}}{|\beta(e)|!}.
\end{equation}

\subsection{First hyperedge splitting of essential walks}
Let us start with necessary definitions. The first vertex $w_1$
of the essential walk $w$ is called the root of the walk. We
denote it by $\rho$. Let us denote the ordered vertices of  $e_1\cap \widetilde{V}_w$ by $e_1(1),e_1(2),\ldots ,e_1(\varkappa)$, where
$\varkappa=|e_1\cap \widetilde{V}_w|$. We denote by $l$ the
walk's length and by $r$ the number of steps of $w$ starting from
root $\rho$.
  In this subsection  we derive the recurrent
  relations by splitting of the walk (or of the tree) into several
  parts. To describe this procedure, it is convenient to consider
    the set of the essential walks of length $l$ such that they
have $r$ steps starting from the root $\rho$.  We denote this set
by $\Lambda(l,r)$. Denote by $S(l,r)$   the sum of contributions of
all the walks of $\Lambda(l,r)$. Let us remove the first hyperedge
$e_1$ from $G_w$ and denote by $G_1,G_2,\ldots, G_\varkappa$ the hypergraphs with roots $e_1(1),e_1(2), \ldots, e_1(\varkappa)$ accordingly
obtained.
 Denote by $u,v\in\mathbb{Z}_+^\varkappa$ the vector, where $u_i$ is  the number of steps along $G_i$ and
  $v_i$ is the number of steps along $G_i$ starting from $e_1(i)$.
Also let's   denote by $f\in\mathbb{Z}_+^\varkappa$ the vector, where $f_i$ is  the number of steps within hyperedge $e_1$  starting from $e_1(i)$. It is clear that the following inequalities hold
for all essential walks (excepting the walk of length zero) $1\leq
f_1\leq r$, $||f||_1+||u||_1= l$. Let us denote by $\Lambda^{(\varkappa)}_2(l,r,f, u,v)$ the
set of the essential walks with fixed parameters $l$, $r$,$\varkappa$,$f$, $u$,
$v$. And  let us denote by $S^{(\varkappa)}_2(l,r,f, u,v)$  the sum of contributions of the walks
of $\Lambda^{(\varkappa)}_2(l,r,f, u,v)$. Denote by $\Lambda_1(l,f)$ the set of
the essential walks of  such that their skeleton
has only one hyperedge and fixed parameter $f$.
  Also we denote by $S_1(l,f)$ the sum of weights   of all   the
walks of $\Lambda_1(l,f)$. Now we can formulate the lemma of
decomposition. It allows express $S$ as  function of  $S$ and $S_2$. Also it allows express $S_2$ as  function of  $K$. And It allows express $K$ as  function of  $K$.

\begin{lemma}[First splitting lemma] The following relations hold  \label{l1}
\begin{equation}\label{l11}
 \mathrm{S}(l,r)\!=\!\sum_{\varkappa=1}^{q}\sum_{\substack{f\in\mathbb{Z}_+^\varkappa\\||f||_1\le l\\f_1\le r\\ \forall i\in [\varkappa]f_i\ge 1 }}
 \sum_{\substack{u\in\mathbb{Z}_+^\varkappa\\||u||_1= l-||f||_1 }}\sum_{\substack{v\in\mathbb{Z}_+^\varkappa\\v\le u\\v_1=r-f_1 }}\mathrm{S}^{(\varkappa)}_2(l,r,f, u,v),
\end{equation}
$$
\mathrm{S}^{(\varkappa)}_2(l,r,f, u,v) \!=\!\frac{p\cdot X_{||f||_1}}{(q-\varkappa)!}\cdot K^{(\varkappa,q)}_{1;f}\cdot {r-1\choose
f_1-1}\cdot  \mathrm{S}(u_1,v_1) \cdot
$$
\begin{equation}\label{l12a}
 \cdot \prod_{i=2}^{\varkappa}{f_i+v_i-1\choose
f_i-1}\cdot  \mathrm{S}(u_i,v_i).
\end{equation}

\end{lemma}
\

\begin{proof}
 The first equality is obvious. The last  equality follows
from the bijection $F$
$$\mathrm{\Lambda}^{(\varkappa)}_2(l,r,f, u,v) \stackrel{bij}{\to}
\mathrm{\Lambda}_{1}(||f||_1,f) \times
  \prod_{i=1}^\varkappa\mathrm{\Lambda}(u_i,v_i) \times $$
\begin{equation}\label{bij}
    \times \prod_{i=1}^\varkappa\mathrm{\Theta}_i(f_i+v_i,f_i) ,
\end{equation}
where $\mathrm{\Theta}_i(f_i+v_i-1,f_i-1)$ is the set of sequences of 0 and 1
of
  length $f_i+v_i-1$ such that there are exactly $f_i-1$ symbols 1 in the
  sequence and the first symbol is 1 if $i=1$ (and the last symbol is 1 if $i\neq 1$).

     Let us construct this
  mapping $F$. Regarding one particular essential walk $w$ of
  $\Lambda^{(\varkappa)}_2(l,r,f, u,v)$, we consider the first hyperedge $e_1$ of the
  hypergraph $G_w$ and separate $w$ in $q+1$ parts, a walk $w^{(0)}$ within $e_1$ and walks $\{w^{(0)}\}_{i=0}^\varkappa$ within hypergraphs $\{G_i(w)\}_{i=1}^\varkappa$ accordingly. We start from root $e_1(1)$. Then we add  special
  codes $\{\theta_i\}_{i=1}^\varkappa$ that determines the transitions at vertices $\{e_1(i)\}_{i=1}^\varkappa$.
  Obviously these $q+1$ parts are walks, but not necessary minimal
  walks. Then  we minimize these walks.
  This decomposition is
  constructed by the following algorithm. We run over $w$ and
  simultaneously draw $\varkappa+1$ parts and codes. If
  the current step belongs to some $G_l$, we add it to the walk $w^{(l)}$. Since $G_w$ is a hypertree, any pair of distinct hypergraphs from $\{G_i\}_{i=0}^\varkappa$ have no common hyperedges. Also
   any pair of distinct hypergraphs from $\{G_i\}_{i=1}^\varkappa$ have no common vertices. And $V_{G_0}\bigcap V_{G_i}=e_1(i)$. These facts ensure the integrity of the walks $\{w^{(i)}\}_{i=1}^\varkappa$. The codes are
  constructed as follows.  Each time the walk leaves the vertex $e_1(l)$ the
  sequence $\theta_l$ is enlarged by one symbol. If
  current step is inside $e_1$ this symbol is "1", otherwise  this symbol is "0".
   It is clear that
  the first element of the sequence $\theta_1$ is "1" (because of definition of the first hyperedge $e_1$ the first step $(w_1=e_1(1),e_1,w_2=e_1(2))$ is inside the first hyperedge $e_1$), the number of signs "1" of the sequence $\theta_1$ is equal to
   $f_1$, and the  full length of the sequence is $r=f_1+v_1$.  It is also clear that for any number $i\in [\varkappa] \setminus {1}$
  the last element of the sequence $\theta_i$ is "1" (because the walk $w$ is closed and $G_w$ is a hypertree), the number of signs "1" of the sequence $\theta_i$ is equal to
   $f_i$, and the  full length of the sequence is $f_i+v_i$. Now we
minimize $\varkappa+1$ parts. Thus, we have constructed
the decomposition of the essential walk $w$ and the mapping $F$.


\begin{example}For
$w=(1,e_1,2, e_3,3,e_3,4,e_3,5,e_3,3,e_2,2,e_4,6,e_4,2,e_1,7, e_5,8,e_5,7,e_6,9,e_6,7,$ $e_1,1,e_7,10,e_7,11,e_7,1,e_1,7,e_6,9,e_6,7,e_5,8,e_5,7,e_1,2,e_1,1,e_1,2,e_1,7,e_1,1,e_7,11,e_7,1,$ $e_1,2,$ $e_1,1)$
the unminimized parts look like
$w^{(0)}=(1,e_1,2,e_1,7,e_1,1,e_1,7,e_1,2,e_1,1,e_1,2,e_1,7,$ $e_1,1,e_1,2,e_1,1)$,
$w^{{(1)}}=(1,e_7,10,e_7,11,e_7,1,e_7,11,e_7,1)$, $w^{{(2)}}=(2,e_2,3,e_3,4,e_3,5,e_3,3,$ $e_2,2,e_4,6,e_4,2)$, $w^{{(3)}}=(7,e_5,8,e_5,7,e_6,9,e_6,7,e_6,9,e_6,7,e_5,8,e_5,7)$,
the codes look like $\theta_1=(1,0,1,1,0,1)$, $\theta_2=(0,0,1,1,1,1)$, $\theta_3=(0,0,1,0,0,1,1)$. And $f=(4,4,3)$,  $v=(2,2,4)$, $u=(4,4,3)$.
\end{example}


 The weight $\theta(w)$ of the
essential walk is  multiplicative with respect to edges.  So $\theta(w)=\prod_{i=1}^\varkappa \theta(w^{(i)})$.  Then the
weight of the essential walk $w$ is equal to the product of
weights of left and right parts. The walk of zero length has unit
weight. The weight of each walk from the set $\mathrm{\Lambda}_{1}(||f||_1,f)$  is $\frac{p\cdot X_{||f||_1}}{(q-\varkappa)!} $. the cardinality of
 the set $\mathrm{\Lambda}_{1}(||f||_1,f)$ equals to $ K^{(\varkappa,q)}_{1;f}$ by definition of the numbers $K$. Combining all these facts with (\ref{bij}), we obtain
\begin{equation}\label{aux}
 \mathrm{S}^{(\varkappa)}_2(l,r,f, u,v)  \!=\!\frac{p\cdot X_{||f||_1}}{(q-\varkappa)!}\cdot K^{(\varkappa,q)}_{1;f}\cdot \prod_{i=1}^\varkappa \left(|\mathrm{\Theta}_i(f_i+v_i,f_i)|\cdot  \mathrm{S}(u_i,v_i)\right) .
\end{equation}
Taking into account that $|\mathrm{\Theta}_{1}(v_i+f_i,f_i)|={v_i+f_i-1\choose
f_i-1}$, we derive (\ref{l12a}) from (\ref{aux}).

By definition, $ K^{(\varkappa,q)}_{j;f}$  is the number of minimal walks along an edge $[\varkappa]$ starting at the root 1, ending at the vertex $j$, passing through all vertices of $[\varkappa]$, and having the edge vector $f$. Let us remove the last step of the minimal walks from the set of walks $ K^{(\varkappa,q)}_{j;f}$. Let us call the resulting set $Q$. These two sets have the same cardinality. The minimal walk from $Q$ can pass all the vertices of $[\varkappa]$ or pass all but $\varkappa$. In the latter case, the last step of the corresponding t walk from $ K^{(\varkappa,q)}_{j;f}$ have to end at $\varkappa$. The edge vector of the minimal walk is   decreased by one in the corresponding row. Thus the following equality $Q=$ is true. Therefore, equalities (\ref{sys_K})-(\ref{ini_cond_K}) hold.

Now let us prove that for any given  elements $\{w^{(i)}\}_{i=0}^\varkappa$ of
$\{\mathrm{\Lambda}(u_i,v_i)\}_{i=1}^\varkappa$, $w^{(0)}$ of
  $\mathrm{\Lambda}_1(||f||_1,f)$, and the sequences $\{\theta_i\in \mathrm{\Theta}_{i}(v_i+f_i,f_i)\}_{i=1}^\varkappa$, one can construct one and only one
  element $w$ of $\mathrm{\Lambda}^{(\varkappa)}_2(l,r,f, u,v)$.
   We do this
  with the following gathering algorithm. First, we take representatives of the class of walks $\{w^{(i)}\}_{i=1}^\varkappa$ in such a way that the pairwise intersections of the sets are empty. Next, we take a representative of the class of walks $w^{(0)}$ in such a way that the property $V_{G_0}\bigcap V_{G_i}=e_1(i)$ is preserved. We start from the root $e_1(1)$. As soon as we get to the vertex $e_1(i)$, we look at the next code element. If it is equal to one, we go through the next step of the first edge walk $w^{(0)}$ and add it to the walk $w$, otherwise we go through the next step of the walk $w^{(i)}$ and add it to the walk. If the next step does not start from any vertex of the first edge, then we are in the  walk $w^{(i)}$, we just take the next step along it and add it to the walk $w$. At the end, we minimize the resulting walk $w$.
   To illustrate
the gathering procedures we give the following example. This example has the same walks $\{w^{(i)}\}_{i=0}^\varkappa$ as the previous example, but different codes $\{\theta_i\}_{i=1}^\varkappa$.


\begin{example} For the unminimized parts
$w^{(0)}=(1,e_1,2,e_1,7,e_1,1,e_1,7,e_1,2,e_1,1,e_1,2,e_1,7,$ $e_1,1,e_1,2,e_1,1)$,
$w^{{(1)}}=(1,e_7,10,e_7,11,e_7,1,e_7,11,e_7,1)$, $w^{{(2)}}=(2,e_2,3,e_3,4,e_3,5,e_3,3,$ $e_2,2,e_4,6,e_4,2)$, $w^{{(3)}}=(7,e_5,8,e_5,7,e_6,9,e_6,7,e_6,9,e_6,7,e_5,8,e_5,7)$,
 and the codes  $\theta_1=(1,1,1,1,0,0)$, $\theta_2=(1,1,1,0,0,1)$, $\theta_3=(1,1,0,0,0,0,1)$ the gathering
procedure gives the unminimized walk  $w=(1,e_1,2,e_1,7,e_1,1,e_1,7, e_1,2,e_1,1,e_1,2,e_1,7,e_5,8,e_5,7,e_6,9,$ $e_6,7,$ $e_6,9,e_6,7,e_5,8,e_5,7,e_1,1,e_1,2,e_2,3,e_3,4,e_3,5,e_3,3,e_2,2,$ $e_4,6,e_4,2,e_1,1,e_7,10,e_7,11,$ $e_7,1,$ $e_7,11,e_7,1)$.
\end{example}

It is clear that the splitting and gathering are  injective
mappings. Their domains are finite sets, and therefore the
corresponding mapping (\ref{bij}) is bijective.

 This completes the
proof of Lemma \ref{l1}. \end{proof}

\subsection{$C^{(N,p,q)}_{k,m}.$}

\qquad Let us denote the set of double closed walks of $k$ and $m$ steps over
the set $[N]$ by
$\mathfrak{D}^{(N)}_{k,m}\stackrel{\rm def}{\equiv}W^{(N)}_{k} \times
W^{(N)}_{m}$.
 For
$\mathfrak{d}\!=\!(w^{(1)},w^{(2)}) \!\in\! \mathfrak{D}^{(N)}_{k,m}$ let us denote
$
a(\mathfrak{d})\!=\!a(w^{(1)})\cdot a(w^{(2)})$,
$d^{(N,p,q)}(\mathfrak{d})\!=\!d^{(N,p,q)}(w^{(1)})\cdot d^{(N,p,q)}(w^{(2)})$.

 Then we obtain
$$
C^{(N,p,q)}_{k,m}\!=\!\frac{1}{N^2} \sum_{\mathfrak{d}=(w^{(1)},w^{(2)})\in
\mathfrak{D}^{(N)}_{k,m}} \left\{ {\mathbb E} a(\mathfrak{}d) \cdot {\mathbb E}
d^{(N,p,q)}(\mathfrak{d}) -\right.
$$
\begin{equation}\label{eq:cor2}
 \left.-{\mathbb E} a(w^{(1)}) \cdot {\mathbb E}
d^{(N,p,q)}(w^{(1)}) \cdot {\mathbb E} a(w^{(2)}) \cdot {\mathbb E}
d^{(N,p,q)}(w^{(2)})\right\}.
\end{equation}

For closed double walks $\mathfrak{d}\!=\!(w^{(1)},w^{(2)}) \!\in\!\mathfrak{D}^{(N)}_{k,m}$
let us denote
$
n_{\mathfrak{d}}(e)=n_{w^{(1)}}(e)+n_{w^{(2)}}(e).
$
Also let us introduce simple  hypergraph $G_{\mathfrak{d}}\!=\!G_{w^{(1)}}\ \cup G_{w^{(2)}}$ for double walk $\mathfrak{d}=(w^{(1)},w^{(2)})\!\in\!\mathfrak{D}^{(N)}_{k,m}$, i.e.
$V_{\mathfrak{d}}\!=\!V_{w^{(1)}}\ \cup V_{w^{(2)}}$ and
$E_{\mathfrak{d}}\!=\!E_{w^{(1)}}\ \cup E_{w^{(2)}}$.
Then, we can rewrite \ref{eq:cor2}  in the following form
$$
C^{(N,p,q)}_{k,m}\!=\! \frac{1}{N^2} \sum_{\mathfrak{d}=(w^{(1)},w^{(2)})\in
\mathfrak{D}^{(N)}_{k,m}} \left\{ \prod_{e\in E_{\mathfrak{dw}}}{\mathbb E}
a_e^{n_{\mathfrak{d}}(e)} \cdot {\mathbb E} \left[d^{(N,p,q)}_e \right]^{n_{\mathfrak{d}}(e)} -\right.
$$
$$ -\left.\prod_{e\in
E_{w{(1)}}}{\mathbb E} a_e^{n_{w^{(1)}}(e)} \cdot {\mathbb E}
\left[ d^{(N,p,q)}_e\right]^{n_{w^{(1)}}(e)}  \cdot \prod_{e\in E_{w{(2)}}}{\mathbb E}
a_e^{n_{w^{(2)}}(e)} \cdot {\mathbb E} \left[ d^{(N,p,q)}_e\right]^{n_{w^{(2)}}(e)}\right\}\!=
$$
$$
=\! \frac{1}{N^2} \sum_{\mathfrak{d}=(w^{(1)},w^{(2)})\in
\mathfrak{D}^{(N)}_{k,m}} \left\{
\left(\frac{p}{N^{q-1}}\right)^{|E_{d}|}\cdot \prod_{e\in E_{d}}X_{n_{\mathfrak{d}}(e)}
 -\right.
$$
\begin{equation}\label{eq:walks}
\left. \left(\frac{p}{N^{q-1}}\right)^{|E_{w^{(1)}}|+|E_{w^{(2)}}|}\cdot
\prod_{e\in E_{w^{(1)}}}X_{n_{w^{(1)}}(e)}
\prod_{e\in E_{w^{(2)}}}X_{n_{w^{(2)}}(e)}\right\}.
\end{equation}

To perform the limiting
transition
   $N\to\infty$ it is natural to separate $\mathfrak{D}^{(N)}_{k,m}$ into classes of
equivalence.
   Double walks  $\mathfrak{d}=(w^{(1)},w^{(2)})$ and $\mathfrak{u}=(u^{(1)},u^{(2)})$ from  $\mathfrak{D}^{(N)}_{k,m}$ are equivalent
 if and only if their first walks  are   equivalent and their second walks  are   equivalent:
 $$
 \ \mathfrak{d} \sim \mathfrak{u} \Leftrightarrow \left(w^{(1)} \sim u^{(1)} \wedge  w^{(2)} \sim u^{(2)}\right).
 $$
 Let us denote by $[\mathfrak{d}]$ the class of equivalence of double walk $\mathfrak{d}$ and by
$\mathfrak{C}^{(N)}_{k,m}$ the set of such classes. Then we can rewrite (\ref{eq:walks})  in the following form
$$
 C^{(N,p,q)}_{k,m}=\!\frac{1}{N^2} \sum_{[\mathfrak{d}]\in \mathfrak{D}^{(N)}_{k,m}}
\left\{\frac{p^{|E_{\mathfrak{d}}|}}{N^{(q-1)|E_{\mathfrak{d}}|}}\cdot\prod_{i=1}^{\kappa} N
\cdot (N-1) \cdot \ldots \cdot (N-|V_{\mathfrak{d}}|+1)\;\cdot \right.
 $$
 $$
 \left.\cdot\left(
\prod_{e\in E_{\mathfrak{d}}}X_{n_{\mathfrak{d}}(e)}\cdot \prodd_{\substack{\alpha\subset 2^{[E_\mathfrak{d}]}\\ \alpha \neq \emptyset}}\frac{1}{|\beta_\alpha(\mathfrak{d})|!} -\frac{p^{|E_{w^{(1)}}|+|E_{w^{(2)}}|-|E_{\mathfrak{d}}|}}
{N^{(q-1)\cdot (|E_{w^{(1)}}|+|E_{w^{(2)}}|-|E_{\mathfrak{d}}|)}} \cdot \prod_{e\in
E_{w^{(1)}}}X_{n_{w^{(1)}}(e)}\cdot \right.\right.
 $$

\begin{equation}\label{eq:class}
\left.\left. \prod_{e\in
E_{w^{(2)}}}X_{n_{w^{(2)}}(e)}\cdot \prodd_{\substack{\alpha\subset 2^{[E_{w^{(1)}}]}\\ \alpha \neq \emptyset}}\frac{1}{|\beta_\alpha(w^{(1)})|!}\cdot \prodd_{\substack{\alpha\subset 2^{[E_{w^{(2)}}]}\\ \alpha \neq \emptyset}}\frac{1}{|\beta_\alpha(w^{(2)})|!}\right)\right\}.
\end{equation}
Let us define a formal order of pass for double walk $\mathfrak{d}=(w^{(1)},w^{(2)})\!\in\!\mathfrak{D}^{(N)}_{k,m}$:
$$
\mathfrak{d}_i=\left\{\begin{array}{ll} w^{(1)}_i, & \textrm{if } 1 \leq
i\leq k \\ w^{(2)}_{i-k}, & \textrm{if } k+1 \leq i\leq k+m.
\end{array} \right.
$$
Let us denote the set of all minimal double walks of $\mathfrak{D}^{(N)}_{k,m}$ by $\mathfrak{M}^{(N)}_{k,m}$.
Then we obtain
  $$
 N\cdot C^{(N,p,q)}_{k,m} \!=\!
\sum_{w\in \mathfrak{M}^{(N)}_{k,m}} \prod_{i=1}^{\kappa}   \lim_{N \to \infty}
\left[\frac{N^{|V_{\mathfrak{d}}|-(q-1)\cdot|E_{\mathfrak{d}}|-1}}
 {p^{-|E_{\mathfrak{d}}|}}\cdot\left(
\prod_{e\in E_{\mathfrak{d}}}X_{n_{\mathfrak{d}}(e)}\cdot \prodd_{\substack{\alpha\subset 2^{[E_\mathfrak{d}]}\\ \alpha \neq \emptyset}}\frac{1}{|\beta_\alpha(\mathfrak{d})|!}  - \right.\right.
$$
\begin{equation}\label{eq:min_lim}
 \left. \left. -\frac{p^{c(\mathfrak{d})}} {N^{c(\mathfrak{d})}} \cdot
\prod_{e\in E_{w^{(1)}}}X_{n_{w^{(1)}}(e)} \prod_{e\in
E_{w^{(2)}}}X_{n_{w^{(2)}}(e)}\cdot \prodd_{\substack{\alpha\subset 2^{[E_{w^{(1)}}]}\\ \alpha \neq \emptyset}}\frac{1}{|\beta_\alpha(w^{(1)})|!}\cdot \prodd_{\substack{\alpha\subset 2^{[E_{w^{(2)}}]}\\ \alpha \neq \emptyset}}\frac{1}{|\beta_\alpha(w^{(2)})|!}\right) \right],
\end{equation}
 where $c(\mathfrak{d})$ is a number of common hyperedges of $G_{w^{(1)}}$ and $G_{w^{(2)}}$, i.e. $c(\mathfrak{d})=|E_{w^{(1)}}|+|E_{w^{(2)}}|-|E_{\mathfrak{d}}|$.

  $\mathfrak{M}_{k,m}$ is a finite set.  $G_{\mathfrak{d}}$ has at most 2 connected components. But if $G_{\mathfrak{d}}$ has exactly 2 connected components then  $V_{w^{(1)}}\cap V_{w^{(2)}}\!=\!\varnothing
\Rightarrow E_{w^{(1)}}\cap E_{w^{(2)}}\!=\!\varnothing
\Rightarrow c(\mathfrak{d})=0$ $ \Rightarrow  \prod_{e\in
E_{\mathfrak{d}}}X_{n_{\mathfrak{d}}(e)}\cdot \prod_{\substack{\alpha\subset 2^{[E_\mathfrak{d}]}\\ \alpha \neq \emptyset}}\frac{1}{|\beta_\alpha(\mathfrak{d})|!} - \frac{p^{c(\mathfrak{d})}} {N^{c(\mathfrak{d})}} \cdot
\prod_{e\in E_{w^{(1)}}}X_{n_{w^{(1)}}(e)} \prod_{e\in
E_{w^{(2)}}}X_{n_{w^{(2)}}(e)}\cdot \prodd_{\substack{\alpha\subset 2^{[E_{w^{(1)}}]}\\ \alpha \neq \emptyset}}\frac{1}{|\beta_\alpha(w^{(1)})|!}\cdot \prodd_{\substack{\alpha\subset 2^{[E_{w^{(2)}}]}\\ \alpha \neq \emptyset}}\frac{1}{|\beta_\alpha(w^{(2)})|!}\!=\! 0.$ Hence the contribution of such minimal double walks to $ N\cdot C^{(N,p,q)}_{k,m}$ equals to 0. Otherwise  $G_{\mathfrak{d}}$ is connected hypergraph, so
 $|V_\mathfrak{d}|-(q-1)\cdot|E_\mathfrak{d}|-1\leq 0$. So (i) of  Theorem \ref{thm:1} is proved.

   \subsection{Carleman's condition.}

The multiplier $C^k$ or $C^{-k}$ does not affect the fulfillment of the Carleman condition (\ref{Carleman}). $m_{2k}^{(p,q)}$ equals the total weight of all essential walks of lenght $2k$.
   Let us define a cluster  of a essential walk $w$ as a set of nonordered hyperedges incident to a given vertex $v\:$ of $G_w$ ($v$ is a center of the cluster). So cluster is a subhypergraph of the skeleton $G_w$. Also let us define an the ordered cluster of a essential walk $w$ as a cluster of the essential walk $w$ with  sequence of numbers $n_w(e_0)$, $n_w(e_1)$,\ldots, $n_w(e_l)$, where $\{e_r\}_{r=0}^{l}$ is a sequence of hyperedges
   of the cluster ordered by time of their first passing. Here we assume that all the numbers are even $n_w(e_0)=2j, n_w(e_1)=2i_1,\ldots,n_w(e_l)=2i_l$.  How to reduce the general case to the even case will become clear later. Let us assign to each cluster a set of numbers that records the number of steps from the center of the cluster $v$ inside the zero hyperedge inside $f_*^{(0)}$, the first hyperedge $f_1^{(1)}$, inside the second hyperedge $f_1^{(2)}$, and so on. It is obvious that $ f_*^{(0)}\le j, f_1^{(1)}\le i_1,\ldots ,f_l^{(l)}\le i_l.$ Number of ways to make an ordered cluster with numbers from an ordered cluster is less than $j\cdot 2^s$. We can derive from the splitting lemma that number of passes covering ordered cluster with numbers $ f_*^{(0)}, f_1^{(1)},\ldots ,f_l^{(l)}$ is less than or equal to $q^{2s}\cdot n(j;i_1,i_2,\ldots,i_l)$, where
  \begin{equation}\label{pass_cluster}
  n(j;i_1,i_2,\ldots,i_l)={j+s-1\choose j-1}\cdot \fracd{s!\prodd_{r=1}^{l}i_r}{s\cdot(s-i_1)\cdot(s-i_1-i_2)\cdot \ldots \cdot i_l\cdot \prodd_{r=1}^{l}i_r!}.
  \end{equation}
  Indeed, each step inside the hyperedge can be chosen in no more than $q$ options. And the number of all choices of code sequences that control transitions between hyperedges for the maximal variant is equal to $n(j;i_1,i_2,\ldots,i_l)$. Indeed, steps to  center of ordered cluster $v$ uniquely determined by choice of steps from the vertex $v$. We can choose steps from $v$ along edge $e_0$ by ${j+s-1\choose j-1}$ ways. After that we can choose steps from $v$ along edge $e_1$ by ${s-1\choose i_1-1}$ ways. After that we can choose steps from $v$ along edge $e_2$ by ${s-i_1-1\choose i_2-1}$ ways and so on. Thus the number of passes covering ordered cluster with numbers $2j$, $2i_1$, $2i_2$, \ldots $2i_l$ ($s=\sum_{j=1}^{l}i_j$)  is equal to
   \begin{equation}\label{pass_cluster1}
  {j+s-1\choose j-1}\cdot {s-1\choose i_1-1}\cdot {s-i_1-1\choose i_2-1}\cdot \ldots \cdot{s-i_1-i_2-\ldots -i_{l-1}-1\choose i_l-1}={j+s-1\choose j-1}\cdot
  \end{equation}
  $$
   \cdot \fracd{s!\cdot i_1}{s\cdot(s-i_1)!\cdot i_1!}\cdot \fracd{(s-i_1)!\cdot i_2}{(s-i_1)\cdot(s-i_1-i_2)!\cdot i_2!}\cdot \ldots\cdot\fracd{(s-i_1-i_2-\ldots-i_{l-1})!\cdot i_l}{(s-i_1-i_2-\ldots-i_{l-1})\cdot 0!\cdot i_l!}
  $$
  But it is evident that numbers in (\ref{pass_cluster}) and (\ref{pass_cluster1}) coincide.

  Let us define the weight of the ordered numbered cluster by
 \begin{equation}\label{weight_cluster}
 \theta(j;i_1,i_2,\ldots,i_l)=q^{2s}\cdot n(j;i_1,i_2,\ldots,i_l)\cdot\prodd_{r=1}^{l} (p\cdot X_{2i_r}).
  \end{equation}
  Combining (\ref{weight_cluster}) with (\ref{Karleman_cond}) and estimate
  $C_1 k^ke^{-k}\le k!\le C_2 k^ke^{-k}$, we obtain the following estimate
  \begin{equation}\label{est_weight_cluster}
 \theta(j;i_1,i_2,\ldots,i_l)\le 2^j q^{2s} C_3^ss^{2s}(1+p)^s.
  \end{equation}
  The estimate \ref{est_weight_cluster} remains valid for the general case. Indeed, If there are two odd numbers in one cluster, we replace them with the nearest even numbers, preserving the total sum. In this case, one of the estimation factors will not change, the other will only increase, and the third will change by no more than a constant factor. If there is only one odd number left in a cluster, then we combine it with an odd number from another cluster with the same effect.

  Let us define an ordered skeleton as a skeleton with all ordered clusters and chosen root. The number of ordered skeletons with $i$ hyperedgees is equal to the corresponding Catalan-Fuss number $\fracd{1}{qi+1}{qn+1 \choose i}$. The number of ways to assign $\{n_w(e_j)\}_{j=1}^i$  to the ordered skeleton is ${2k+i-1\choose i-1}$. So the number of all ordered numbered skeleton is not greater than $\sumd_{i=0}^k\left(\fracd{1}{qi+1}{qn+1 \choose i}\cdot {2k+i-1\choose i-1}\right)\le 2^{(q+3)k}$. Let us define the weight of  an ordered skeleton as the number of passes covering it multiplied by $\prodd_{e\in E_w}(p\cdot X_{n_w(e)})$. From (\ref{est_weight_cluster}) we get that the weight of ordered skeleton of essential walk $w$ from $\mathfrak{E}_{2k}$ is not greater than $2^{(q+4)k}q^{2k} {C_3(p+1)}^k k^{2k}$.  So (iii) of Theorem \ref{thm:1} is proved.

  \bibliographystyle{amcjoucc}
\bibliography{amcexample}

\end{document}